\documentclass[a4paper,11pt, english, parskip=full]{scrartcl}
\usepackage[top=2.5cm, bottom=2.5cm, left=2.5cm, right=2.5cm]{geometry}
\usepackage[utf8x]{inputenc}
\usepackage[T1]{fontenc}
\usepackage[english]{babel}
\usepackage{lmodern}

\usepackage[babel]{csquotes}

\usepackage{amsmath}
\usepackage{amssymb}
\usepackage[amsmath,thmmarks]{ntheorem}
\DeclareMathAlphabet{\mathpzc}{OT1}{pzc}{m}{it}
\usepackage{dsfont}
\usepackage{chngcntr}
\counterwithout{equation}{section}

\let\existsorig\exists 
\renewcommand{\exists}{\ \existsorig\ } 
\let\forallorig\forall 
\renewcommand{\forall}{\ \forallorig\ } 

\mathchardef\ordinarycolon\mathcode`\: 
\mathcode`\:=\string"8000 
\begingroup \catcode`\:=\active  
\gdef:{\mathrel{\mathop\ordinarycolon}} 
\endgroup

\theoremstyle{plain}
\theorembodyfont{\upshape}
\newtheorem{Satz}{Theorem}[section]
\newtheorem{Lemma}[Satz]{Lemma}
\newtheorem{Folgerung}[Satz]{Corollary}

\theorembodyfont{\upshape}
\newtheorem{Definition}[Satz]{Definition}
\newtheorem{Bemerkung}[Satz]{Comment}
\newtheorem{Bemerkungen}[Satz]{Comments}

\theoremstyle{break}

\theorembodyfont{\upshape}

\newtheorem{BemerkungenB}[Satz]{Comments}

\theoremstyle{nonumberplain}
\theorembodyfont{\upshape}
\theoremheaderfont{\scshape}
\theoremsymbol{\ensuremath{\hfill\square}}
\theoremseparator{.}

\newtheorem{Beweis}{Proof}

\newcommand{\IC}{\mathbb{C}}
\newcommand{\IN}{\mathbb{N}}
\newcommand{\IR}{\mathbb{R}}
\newcommand{\IZ}{\mathbb{Z}}
\newcommand{\IQ}{\mathbb{Q}}
\newcommand{\IP}{\mathbb{P}}
\newcommand{\dd}{\mathrm{d}}
\newcommand{\el}{\in}
\newcommand{\ra}{\rightarrow}

\newcommand{\intt}{\int\limits}
\newcommand{\summe}{\sum\limits}
\newcommand{\pro}{\prod\limits}
\newcommand{\limes}{\lim\limits}
\newcommand{\e}[1]{{\textbf e}\left(#1\right)}
\newcommand{\ggt}[2]{\left(#1;#2\right)}
\newcommand{\tf}{{\bf d}}
\newcommand{\abs}[1]{\left|#1\right|}
\newcommand{\norm}[1]{\left|\left|#1\right|\right|}
\newcommand{\zweinorm}[1]{\left|\left|#1\right|\right|_2}
\newcommand{\teiltnicht}{\nmid}
\newcommand{\ssumme}[2]{\summe\limits_{#1}\!{\vphantom{\sum}}^*{#2}}
\newcommand{\floor}[1]{\left\lfloor#1\right\rfloor}
\newcommand{\iverson}[1]{\left[#1\right]}
\newcommand{\modu}[3]{{#1\equiv#2~\left(#3\right)}}
\newcommand{\moduu}[3]{{#1\equiv#2\left(#3\right)}}
\newcommand{\nmodu}[3]{{#1\not\equiv#2~\left(#3\right)}}
\newcommand{\nmoduu}[3]{{#1\not\equiv#2\left(#3\right)}}
\newcommand{\Ezwei}[2]{\mathcal{E}_{#1}(#2)}
\newcommand{\Ea}[1]{\Ezwei{a}{#1}}
\newcommand{\EakeineKlammern}{\mathcal{E}_a}
\newcommand{\G}{\operatorname{G}}
\renewcommand{\H}{\operatorname{H}}
\newcommand{\E}{\operatorname{E}}
\renewcommand{\S}{\operatorname{S}}
\newcommand{\T}{\operatorname{T}}
\newcommand{\Ss}{\operatorname{S}^*}
\newcommand{\Ts}{\operatorname{T}^*}
\newcommand{\smallo}{\operatorname{o}}

\title{The asymptotic behavior of limit-periodic functions on primes and an application to $k$-free numbers}
\author{Markus Hablizel}
\date{}

\begin{document}

\maketitle

\begin{abstract}
We use the circle method to evaluate the behavior of limit-periodic functions on primes. For those limit-periodic functions that satisfy a kind of Barban-Davenport-Halberstam condition and whose singular series converge fast enough, we can evaluate their average value on primes. As an application, this result is used to show how tuples of different $k$-free numbers behave when linear shifts are applied. 
\end{abstract}

\section{Introduction}\label{Grundlagen_Kapitel}
Limit-periodic functions are those arithmetical functions $f:\IN\ra\IC$ which appear as limits of periodic functions with regards to the \textit{Besicovitch-seminorm} defined via 
$$\norm{f}_2:=\left(\limsup_{x\ra\infty}{\frac{1}{x}\summe_{n\leq x}{\abs{f(n)}^2}}\right)^{\frac{1}{2}}$$
They have some wonderful properties, e.g., the mean-value as well as the mean-value in residue classes exist always. Limit-periodic functions are a special case of almost-periodic functions that have been explored by Harald Bohr und Abram Besicovitch in the 1920s. The distinction between those two lies in the approximation type: Almost-periodic functions appear as limits with regards to the Besicovitch-seminorm of linear combinations of the functions $k\mapsto e^{2\pi i\alpha k}$ with $\alpha\el\IR$, whereas for limit-periodic functions, only $\alpha\el\IQ$ is admissable. For further reference on limit-periodic functions, see \cite{Schwarz_und_Spilker:Arithmetical_Functions}.

The main result of this paper is a statement on the behavior of a limit-periodic function on primes on average. We prove in theorem \ref{Hauptsatz_der_DA} that under certain conditions the asymptotic relation
$$\summe_{p\leq x}{f(p)}=c_f\frac{x}{\log{x}}+\smallo{\left(\frac{x}{\log{x}}\right)}$$
holds, with a constant $c_f$ explicitely given through an infinite series. Brüdern \cite{Bruedern_AdditiveProblems} has considered this result in a more general context.

As an application we show for arbitrary $\alpha_i\el\IN_0$ and $r_i\el\IN_{>1}$
\begin{equation}\label{Gleichung:Aussage_des_Beispiels}
\summe_{p\leq x}{\mu_{r_1}(p+\alpha_1)\cdot\ldots\cdot\mu_{r_s}(p+\alpha_s)}=\pro_{p}{\left(1-\frac{\operatorname{D}^{*}(p)}{\varphi(p^{r_s})}\right)}\,\frac{x}{\log{x}}+{\rm o}\left(\frac{x}{\log{x}}\right)
\end{equation}
where $\mu_k$ denotes the characteristic function of the $k$-free numbers and $\operatorname{D}^{*}(p)$ is a computable function of the prime $p$, depending on the choice of the numbers $\alpha_i$ and $r_i$.

\section{Some basic facts}
We state some basic facts and notation for the later discourse.

\begin{Definition}[$k$-free numbers]\label{DefinitionMueK} For given $k\el\IN_{>1}$ the function $\mu_k$ denotes the characteristic function of the set of \emph{$k$-free numbers}, i.e. 
\begin{equation*}
\mu_k(n):=\begin{cases}0&\mbox{there is a }p\el\IP\mbox{ with }p^k|n\\1&\mbox{otherwise}\end{cases}
\end{equation*}
which is multiplicative. On prime powers it has the values
\begin{equation}\label{MuekaufPrimzPot}
\mu_k(p^r)=1-\iverson{k\leq r}
\end{equation}
where $\iverson{A}$ shall denote the \textit{Iverson bracket} to the statement $A$, i.e., it equals $1$ if $A$ is true, and $0$ otherwise. As it is long known we also have 
\begin{equation}\label{MueKIdentitaet}
\mu_k(n)=\summe_{d^k|n}{\mu(d)}
\end{equation}
\end{Definition}

\begin{Lemma}\label{Lemma:MuQuadrat_durch_Phi_Summe}
For $x\el\IR_{>1}$ we have the asymptotic relation
\begin{equation*}
\summe_{n\leq x}{\frac{\mu(n)^2}{\varphi(n)}}=\log{x}+\gamma+\summe_{p}{\frac{\log{p}}{p\left(p-1\right)}}+{\rm O}\left(\frac{\log{x}}{\sqrt{x}}\right)
\end{equation*}
with the Euler-Mascheroni-constant $\gamma$, see \cite{Bruedern:Einfuehrung_in_die_analytische_Zahlentheorie}.
\end{Lemma}

\begin{Definition}\label{E-FKT}We define the function ${\bf e}$ with period $1$ as usual through
\begin{equation*}
{\bf e}:\IR\ra\left\{z\el\IC:~\abs{z}=1\right\},~x\mapsto e^{2{\rm\pi i}x}
\end{equation*}
We sometimes write ${\bf e}_{\frac{a}{q}}$ for the function $n\mapsto{\bf e}(\frac{an}{q})$. For ${q\el\IN}$ \emph{Ramanujan's sum} ${{\bf c}_q}$ is given by
\begin{equation*}
{\bf c}_q(n):=\ssumme{a\leq q}{{\bf e}_{\frac{a}{q}}(n)}
\end{equation*}
where the star on the sum shall denote the sum over all $a\leq q$ prime to $q$ only, i.e., their greatest common divisor equals $1$.
\end{Definition}

\begin{Lemma}
With the geometric series we have for $a,b\el\IZ$, $0\leq a<b$, $\beta\el\IR$ the inequality
\begin{equation}\label{Gleichung:Exp_von_Beta_Abschaetzung}
\left|\summe_{a<n\leq b}{\e{\beta n}}\right|\leq\min{\left(b-a,\frac{1}{2\norm{\beta}}\right)}
\end{equation}
where $\norm{\beta}$ denotes the distance to the nearest integer. For a proof, see \cite{Kraetzel}.
\end{Lemma}

\subsection*{The space $\mathcal{D}^2$ of limit-periodic functions}\label{Subsection_Grenzperiodische_Funktionen}
For $q\el\IN$, let $\mathcal{D}_q$ be the set of all $q$-periodic functions and $\mathcal{D}:=\bigcup_{q=1}^{\infty}{\mathcal{D}_q}$. Write $\mathcal{D}^2$ for the closure of $\mathcal{D}$ with regards to the Besicovitch-seminorm $\norm{\cdot}_2$ which makes it a normed vector space in a canonical way. Limit-periodic functions are exactly the elements of this vector space.

\begin{Satz}
The vector spaces ${\mathcal{D}_q}$ and $\mathcal{D}$ possess the following bases, see \cite{Schwarz_und_Spilker:Arithmetical_Functions},
\begin{align}
\mathcal{D}_q&=\left<{\bf e}_{\frac{a}{q}}~:~1\leq a\leq q\right>_{\IC}\nonumber\\
\mathcal{D}&=\left<{\bf e}_{\frac{a}{q}}~:~1\leq a\leq q,\,q\el\IN,\,\left(a;q\right)=1\right>_{\IC}\label{Orthonormal_Basis}
\end{align}
\end{Satz}

\begin{Definition}[Besicovitch-seminorm]\label{BesicovitchNorm}
The \emph{Besicovitch-seminorm} of a function ${f\el\mathcal{D}}$ is given through
\begin{equation*}
\zweinorm{f}:=\left(\limsup_{x\ra\infty}{\frac{1}{x}\summe_{n\leq x}{\left|f(n)\right|^2}}\right)^{\frac{1}{2}}
\end{equation*}
Note that for a $q$-periodic function we have the identity
\begin{equation*}
\summe_{n\leq x}{\left|f(n)\right|^2}=\left(\floor{\frac{x}{q}}+{\rm O}\left(1\right)\right)\summe_{n\leq q}{\left|f(n)\right|^2}
\end{equation*}
\end{Definition}

\begin{Bemerkungen}\label{Bemerkungen_Grenzp_von_Shifted_Functions}
If $f$ is a limit-periodic function, so is $\abs{f}$, $\operatorname{Re}(f)$ and $\operatorname{Im}(f)$ as well as with ${a\el\IZ,}$\,${b\el\IN,}$
\begin{equation*}
\begin{aligned}
n&\mapsto f(n+a)\\
n&\mapsto f(bn)
\end{aligned}
\end{equation*}

Furthermore, the \textit{mean-value} 
\begin{equation*}
\operatorname{M}(f):=\limes_{x\ra\infty}{\frac{1}{x}\summe_{n\leq x}{f(n)}}
\end{equation*}
exists for every $f\el\mathcal{D}^2$, as well as the \textit{mean-value in residue classes}
\begin{equation}\label{DefvonEta}
\eta_f(q,b):=\limes_{x\ra\infty}{\frac{1}{x}\summe_{\substack{n\leq x\\\moduu{n}{b}{q}}}{f(n)}}
\end{equation}
 for arbitrary ${b,q\el\IN}$. For the respective proofs, see \cite{Schwarz_und_Spilker:Arithmetical_Functions}.
\end{Bemerkungen}
 
\begin{Bemerkung}\label{Bemerkung:Abschaetzung_fuer_eta}
For ${f\el\mathcal{D}^2}$ we have ${\eta_f(q,b)\ll q^{-\frac{1}{2}}}$ which is easily seen with the Cauchy-Schwarz-inequality
\begin{equation*}
\frac{1}{x}\abs{\summe_{\substack{n\leq x\\\moduu{n}{b}{q}}}{f(n)}}\leq\left(\frac{1}{x}\summe_{n\leq x}{\abs{f(n)}^2}\right)^{\frac{1}{2}}\left(\frac{1}{x}\summe_{\substack{n\leq x\\\moduu{n}{b}{q}}}{1}\right)^{\frac{1}{2}}\ll\left(\frac{1}{x}\left(\frac{x}{q}+1\right)\right)^{\frac{1}{2}}\ll q^{-\frac{1}{2}}
\end{equation*}
where $\ll$ denotes as usual Vinogradov's symbol.
\end{Bemerkung}

\begin{Lemma}[Parseval's identity]\label{Parseval_Lemma}
As the basis (\ref{Orthonormal_Basis}) is an orthonormal basis of~${\mathcal{D}^2}$, Parseval's identity holds as well
\begin{equation*}
\summe_{q=1}^{\infty}{\ssumme{a\leq q}{\left|\operatorname{M}(f\cdot\overline{{\bf e}_{\frac{a}{q}}})\right|^2}}=\zweinorm{f}^2
\end{equation*}
\end{Lemma}

The following example of a limit-periodic function is used in the application at the end of this paper.
\begin{Lemma}\label{muekistgrenzperiodisch}
The function $\mu_k$ is not periodic, but it is limit-periodic.
\begin{Beweis}
Assume we have a natural number $R$ with
$$\mu_k(n+R)=\mu_k(n)$$
for all $n\el\IN$. Then we can deduce that for each $p\el\IP$ and~$m\el\IN$ we have
\begin{equation}
\label{phochkteiltnicht}
p^k\teiltnicht(1+mR)
\end{equation}
which is easily seen to be false with the theorem of Fermat-Euler.
For the proof of the limit-periodic property, define for $k\el\IN$, $k\geq 2$, and $y\el\IR_{>2}$ the arithmetical function ${\varkappa_k^{\left(y\right)}}$ through
\begin{equation*}
\varkappa_k^{\left(y\right)}(n):=\begin{cases}\mu(s)\pro_{p|n}{\iverson{p\leq y}}&\mbox{if }n=s^k\mbox{ with}s\el\IN\\0&\mbox{otherwise}\end{cases}
\end{equation*}
Then ${\varkappa_k^{\left(y\right)}}$ is multiplicative. As a \textit{Dirichlet-convolution} of multiplicative functions, the function
\begin{equation}\label{DefMueKY}
\mu_k^{\left(y\right)}:=\varkappa_k^{\left(y\right)}\,*\mathds{1}
\end{equation}
is multiplicative as well. It is an approximation to $\mu_k$ as can be seen, when evaluated on prime powers:
\begin{equation}\label{MueKYaufPrimzPot}
\mu_k^{\left(y\right)}\left(p^r\right)=1+\summe_{j\leq r}{\begin{Bmatrix}\mu(p)\iverson{p\leq y}&\mbox{ }j=k\\0&\mbox{ otherwise}\end{Bmatrix}}
=1-\iverson{p\leq y}\iverson{k\leq r}
\end{equation}
which means
\begin{equation*}
\mu_k(n)=\summe_{d^k|n}{\mu(d)}~~~~~~~~\mu_k^{\left(y\right)}(n)=\summe_{\substack{d^k|n\\\forallorig p|d\,:\,p\leq y}}{\mu(d)}
\end{equation*}
With the equations~(\ref{MuekaufPrimzPot}) and (\ref{MueKYaufPrimzPot}) we get for all ${p\el\IP}$ and ${r\el\IN_0}$
$$\mu_k^{\left(y\right)}(p^r)\geq\mu_k(p^r)$$
As both functions are multiplicative and can only attain the values $0$ or $1$, we get directly for all ${n\el\IN}$
$$\mu_k^{\left(y\right)}(n)\geq\mu_k(n)$$
and
\begin{equation}\label{NullundEins}
\left(\mu_k^{\left(y\right)}(n)-\mu_k(n)\right)\el\left\{0,1\right\}
\end{equation}

The function ${\mu_k^{\left(y\right)}}$ is periodic with period ${\mathfrak{r}:=\pro_{p\leq y}{p^k}}$, as
\begin{equation*}
\mu_k^{\left(y\right)}(n)=\summe_{\substack{d^k|n\\\forallorig p|d\,:\,p\leq y}}{\mu(d)}=\summe_{\substack{d^k|\left(n+\mathfrak{r}\right)\\\forallorig p|d\,:\,p\leq y}}{\mu(d)}=\mu_k^{\left(y\right)}(n+\mathfrak{r})
\end{equation*}
From (\ref{NullundEins}) we deduce
\begin{equation*}
\zweinorm{\mu_k^{\left(y\right)}-\mu_k}^2=\limsup_{x\ra\infty}{\frac{1}{x}\summe_{n\leq x}{\left(\summe_{\substack{d^k|n\\\existsorig p|d\,:\,p>y}}{\mu(d)}\right)^2}}=\limsup_{x\ra\infty}{\frac{1}{x}\summe_{n\leq x}{\summe_{\substack{d^k|n\\\existsorig p|d\,:\,p>y}}{\mu(d)}}}
\end{equation*}
and it follows
\begin{align*}
\zweinorm{\mu_k^{\left(y\right)}-\mu_k}^2
&\leq\limsup_{x\ra\infty}{\frac{1}{x}\summe_{n\leq x}{\summe_{\substack{d^k|n\\\existsorig p|d\,:\,p>y}}{1}}}=\limsup_{x\ra\infty}{\frac{1}{x}\summe_{\substack{d^k\leq x\\\existsorig p|d\,:\,p>y}}{\summe_{\substack{n\leq x\\d^k|n}}{1}}}\\
&\leq\limsup_{x\ra\infty}{\frac{1}{x}\summe_{\substack{d\leq x\\\existsorig p|d\,:\,p>y}}{\frac{x}{d^k}}}=\limsup_{x\ra\infty}{\summe_{\substack{d\leq x\\\existsorig p|d\,:\,p>y}}{d^{-k}}}\\
&\leq\summe_{d>y}{d^{-k}}~\stackrel{y\ra\infty}{\longrightarrow}~0
\end{align*}
Thereby, we get
\begin{equation}\label{MueK_konvergiert_gegen_Mue}
\zweinorm{\mu_k-\mu_k^{\left(y\right)}}~\stackrel{y\ra\infty}{\longrightarrow}~0
\end{equation}
which shows the limit-periodic property of $\mu_k$.
\end{Beweis}
\end{Lemma}

\begin{Lemma}\label{Lemma:Grenzp_und_beschr_Fkt}
For a limit-periodic function $f$ that is bounded in addition, the function ${\mu_k f}$ is limit-periodic as well, which can be easily seen.
\end{Lemma}

\subsection*{The Barban--Davenport--Halberstam theorem for $\mathcal{D}^2$}
The Barban-Davenport-Halberstam theorem in its original form for primes proves that the error term in the prime number theorem for arithmetic progressions is small in the quadratic mean, see \cite{Bruedern:Einfuehrung_in_die_analytische_Zahlentheorie}, and for further references \cite{Hooley_3}, \cite{Hooley_9}. We need a corresponding version for limit-periodic functions.

Define the error term in the sum over arithmetic progressions via
\begin{equation}\label{Def_von_E}
\operatorname{E}_f(x;q,b):=\summe_{\substack{n\leq x\\\moduu{n}{b}{q}}}{f(n)}-x\eta_f(q,b)
\end{equation}
Then the following lemma due to Hooley \cite{Hooley_10} holds.
\begin{Lemma}\label{Satz:BDH_fuer_grenzperiodische_Funktionen}
If for all $A\el\IR$, $b,q\el\IN$
$$\operatorname{E}_f(x;q,b)\ll\frac{x}{\left(\log{x}\right)^A}$$
where the implicit constant in Vinogradov's symbol is at most dependent on $A$ or $f$, then we have for all $A\el\IR$ and $Q\el\IR_{>0}$
$$\summe_{q\leq Q}{\summe_{b\leq x}{\abs{\operatorname{E}_f(x;q,b)}^2}}\ll Qx+\frac{x^2}{\left(\log{x}\right)^A}$$
\end{Lemma}
\clearpage
\section{Proof of the main theorem with the circle method}\label{Main_Chapter}
In this section we state and prove the main theorem \ref{Hauptsatz_der_DA} with the circle method of Hardy and Littlewood \cite{Vaughan}. Let $f\el\mathcal{D}^2$ be throughout this section a given function.
\begin{Definition}
We define for $a,q\el\IN$ the \emph{Gaußian sum} of $f$ via
\begin{equation*}
\G_f(q,a):=\summe_{b\leq q}{\eta_f(q,b)\,\e{\frac{ab}{q}}}
\end{equation*}
\end{Definition}

\begin{Bemerkungen} The Gaußian sum is a mean-value, as
\begin{align*}
\G_f(q,a)&=\summe_{b\leq q}{\limes_{x\ra\infty}{\frac{1}{x}\summe_{\substack{n\leq x\\\moduu{n}{b}{q}}}{f(n)}}\,\e{\frac{ab}{q}}}=\limes_{x\ra\infty}{\frac{1}{x}\summe_{b\leq q}{\summe_{\substack{n\leq x\\\moduu{n}{b}{q}}}{f(n)\,\e{\frac{an}{q}}}}}\\
&=\limes_{x\ra\infty}{\frac{1}{x}\summe_{n\leq x}{f(n)\,\e{\frac{an}{q}}}}=\operatorname{M}(f\cdot{\bf e}_{\frac{a}{q}})
\end{align*}
With Parseval's identity we also have
\begin{equation}
\begin{aligned}\label{Parseval_mit_G_f}
\mathfrak{S}_f:=&\summe_{q=1}^{\infty}{\ssumme{a\leq q}{\left|\G_f(q,a)\right|^2}}=\summe_{q=1}^{\infty}{\ssumme{a\leq q}{\left|\summe_{b\leq q}{\eta_f(q,b)\,\e{\frac{ab}{q}}}\right|^2}}\\
=&\summe_{q=1}^{\infty}{\ssumme{a\leq q}{\left|\summe_{b\leq q}{\eta_f(q,b)\,\e{-\frac{ab}{q}}}\right|^2}}\\
=&\summe_{q=1}^{\infty}{\ssumme{a\leq q}{\left|\operatorname{M}(f\cdot\overline{{\bf e}_{\frac{a}{q}}})\right|^2}}=\zweinorm{f}^2
\end{aligned}
\end{equation}
Therefore, for limit-periodic functions the identity $\mathfrak{S}_f=\zweinorm{f}^2$ holds.
The series $\mathfrak{S}_f$ is called \textit{singular series} of $f$.
\end{Bemerkungen}

We are now able to state the main theorem of this paper.
\begin{Satz}\label{Hauptsatz_der_DA}
Let $f\el\mathcal{D}^2$ be an arithmetical function with
\begin{equation}\label{ERSTE_ANFORDERUNG}
\summe_{n\leq x}{\abs{f(n)}^2}=x\zweinorm{f}^2+{\rm o}\left(\frac{x}{\log{x}}\right)
\end{equation}
and the remainder of the corresponding singular series (\ref{Parseval_mit_G_f}) satisfies
\begin{equation}\label{Gleichung:Forderung_an_Betragsquadrat_ueber_Gausssche_Summe}
\summe_{q>\operatorname{w}}{\ssumme{a\leq q}{\left|\G_f(q,a)\right|^2}}={\rm o}\left(\operatorname{w}^{-\frac{1}{r}}\right)
\end{equation}
with $r\el\IR_{>1}$. We then set $Q=Q(x):=\left(\log{x}\right)^{r}$. Furthermore, we demand for all $A\el\IR$
\begin{equation}\label{ZWEITE_ANFORDERUNG}
\summe_{q\leq Q}{\max_{1\leq b\leq q}{\abs{\E_f(x;q,b)}}}\ll\frac{x}{\left(\log{x}\right)^A}
\end{equation}
where the implicit constant in Vinogradov's symbol is at most dependent on $A$.

Then we have
\begin{equation}\label{Gleichung:Haupt_Aussage}
\summe_{p\leq x}{f(p)}=c_f\,\frac{x}{\log{x}}+{\rm o}\left(\frac{x}{\log{x}}\right)
\end{equation}
with a constant $c_f$ that is represented through the infinite series
\begin{equation*}
c_f:=\summe_{q=1}^{\infty}{\frac{\mu(q)}{\varphi(q)}\ssumme{a\leq q}{\G_f(q,a)}}
\end{equation*}
\end{Satz}

\begin{BemerkungenB}
\begin{enumerate}
\item{The condition (\ref{ZWEITE_ANFORDERUNG}) implies
\begin{equation*}
\E_f(x;q,b)\ll\frac{x}{\left(\log{x}\right)^A}
\end{equation*}
for all $A\el\IR,\,b,q\el\IN$, and we can apply theorem \ref{Satz:BDH_fuer_grenzperiodische_Funktionen}.}
\item{The following identity can be verified easily
\begin{equation*}
c_f=\summe_{q=1}^{\infty}{\frac{\mu(q)}{\varphi(q)}\operatorname{M}(f\cdot{\bf c}_q)}
\end{equation*}}
\end{enumerate}
\end{BemerkungenB}
In what follows, we assume the conditions of theorem~\ref{Hauptsatz_der_DA}. For notational simplification, we write $\G,~\E,~\eta$, etc. instead of $\G_f,~\E_f,~\eta_f$, etc. 
\subsection{Split in major and minor arcs}
\begin{Definition}[Major and minor arcs]
With the unit interval $\mathcal{U}:=\left(\frac{Q}{x},1+\frac{Q}{x}\right]$ we define for $a,q\el\IN$ with $1\leq a\leq q\leq Q$, $\left(a;q\right)=1$ the \emph{major arcs} through
\begin{equation*}
\mathfrak{M}(q,a):=\left\{\alpha\el\mathcal{U}~:~\abs{\alpha-\frac{a}{q}}\leq\frac{Q}{x}\right\}
\end{equation*}
Let the symbol $\mathfrak{M}$ denote the union of all major arcs
\begin{equation*}
\mathfrak{M}:=\bigcup_{q\leq Q}{\bigcup_{\substack{a\leq q\\\left(a;q\right)=1}}^{}{\mathfrak{M}(q,a)}}
\end{equation*}
We define the \emph{minor arcs} as usual as the complement in the unit interval
\begin{equation*}
\mathfrak{m}:=\mathcal{U}\setminus\mathfrak{M}
\end{equation*}
\end{Definition}
For sufficient large $x$ each pair of major arcs is disjunct.

\begin{Definition}
We define exponential sums $\S$ and $\T$ for $\alpha\el\IR$ via
\begin{align*}
\S(\alpha)&:=\summe_{n\leq x}{f(n)\,\e{\alpha n}}\\
\T(\alpha)&:=\summe_{p\leq x}{\e{-\alpha p}}
\end{align*}
\end{Definition}
and have then
\begin{equation}
\begin{aligned}\label{Gleichung:Hauptaussage_Zerlegung}
\summe_{p\leq x}{f(p)}&=\intt_{0}^{1}{\S(\alpha)\T(\alpha)\,\dd\alpha}\\
&=\intt_{\mathcal{U}}{\S(\alpha)\T(\alpha)\,\dd\alpha}=\intt_{\mathfrak{M}}{\S(\alpha)\T(\alpha)\,\dd\alpha}+\intt_{\mathfrak{m}}{\S(\alpha)\T(\alpha)\,\dd\alpha}
\end{aligned}
\end{equation}

\subsection{The major arcs}
On~$\mathfrak{M}$ we approximate $\S$ resp. $\T$ by the functions $\Ss$ resp. $\Ts$ that are defined for $\alpha\el\mathfrak{M}(q,a)$, $\alpha=\frac{a}{q}+\beta$, via
\begin{align*}
\Ss(\alpha)&:=\G(q,a)\summe_{n\leq x}{\e{\beta n}}\\
\Ts(\alpha)&:=\frac{\mu(q)}{\varphi(q)}\summe_{2\leq n\leq x}{\frac{\e{-\beta n}}{\log{n}}}
\end{align*}

\begin{Lemma}\label{Lemma:Ts_Abschaetzung}
The function $\Ts$ satisfies on ${\mathfrak{M}(q,a)}$ with $\alpha=\frac{a}{q}+\beta$ the inequality
\begin{equation}\label{Gleichung:Ts_Abschaetzung}
\Ts(\alpha)\ll\frac{\mu(q)^2}{\varphi(q)\log{x}}~\frac{x}{1+\norm{\beta}x}
\end{equation}
\begin{Beweis}
The case $\beta\el\IZ$ is trivial. For $\beta\notin\IZ$ the method of partial summation and estimate (\ref{Gleichung:Exp_von_Beta_Abschaetzung}) can be applied:
\begin{align*}
\Ts(\alpha)&\ll\frac{\mu(q)^2}{\varphi(q)}\left(\frac{1}{\log{x}}\abs{\summe_{2\leq n\leq x}{\e{-\beta n}}}+\intt_{t=2}^{x}{\frac{1}{t\left(\log{t}\right)^2}\abs{\summe_{2\leq n\leq t}{\e{-\beta n}}}\dd t}\right)\\
&\ll\frac{\mu(q)^2}{\varphi(q)}\left(\frac{1}{\norm{\beta}\log{x}}+\frac{1}{\norm{\beta}}\intt_{t=2}^{x}{\frac{1}{t\left(\log{t}\right)^2}\dd t}\right)\ll\frac{\mu(q)^2}{\varphi(q)}\frac{1}{\norm{\beta}\log{x}}
\end{align*}
\end{Beweis}
\end{Lemma}

The next lemma makes the approximation through $\Ss$ and $\Ts$ on the major arcs more precise.
\begin{Lemma}\label{Lemma:S_Gleich_S_Stern_Plus_Rest} We have for $\alpha\el\mathfrak{M}$, $\alpha=\frac{a}{q}+\beta$ and arbitrary $A\el\IR$
\begin{align*}
\S(\alpha)&=\Ss(\alpha)+\summe_{b\leq q}{\e{\frac{ab}{q}}\Xi\left(x;q,b;\beta\right)}\\
\T(\alpha)&=\Ts(\alpha)+{\rm O}\left(\frac{x}{\left(\log{x}\right)^A}\right)
\end{align*}
where
\begin{equation}\label{Def_von_Xi}
\Xi\left(x;q,b;\beta\right):=\e{\beta x}\E(\floor{x};q,b)-2{\rm\pi i}\beta\intt_{t=1}^{x}{\e{\beta t}\E(\floor{t};q,b)\,\dd t}
\end{equation}
\begin{Beweis}
If we evaluate $\S$ at the rational number $\frac{a}{q}$ we get with definition (\ref{Def_von_E})
\begin{align*}
\summe_{n\leq x}{f(n)\,\e{\frac{an}{q}}}&=\summe_{b\leq q}{\e{\frac{ab}{q}}\summe_{\substack{n\leq\floor{x}\\\moduu{n}{b}{q}}}{f(n)}}=\summe_{b\leq q}{\e{\frac{ab}{q}}\left(\floor{x}\eta(q,b)+\E(\floor{x};q,b)\right)}\\
&=\floor{x}\G(q,a)+\summe_{b\leq q}{\e{\frac{ab}{q}}\E(\floor{x};q,b)}
\end{align*}
Applying partial summation twice yields the stated claim:
\begin{align*}
\S\left(\alpha\right)&=\e{\beta x}\S\left(\frac{a}{q}\right)-2{\rm\pi i}\beta\intt_{t=1}^{x}{\e{\beta t}\left(\floor{t}\G(q,a)+\summe_{b\leq q}{\e{\frac{ab}{q}}\E(\floor{t};q,b)}\right)\dd t}\\
&=\G(q,a)\left(\floor{x}\e{\beta x}-2{\rm\pi i}\beta\intt_{t=1}^{x}{\floor{t}\e{\beta t}\,\dd t}\right)+\summe_{b\leq q}{\e{\frac{ab}{q}}\Xi\left(x;q,b;\beta\right)}\\
&=\Ss(\alpha)+\summe_{b\leq q}{\e{\frac{ab}{q}}\Xi\left(x;q,b;\beta\right)}
\end{align*}
For the second statement, we use partial summation another time
\begin{align*}
\T(\alpha)&=\frac{1}{\log{x}}\summe_{p\leq x}{\e{-\alpha p}\log{p}}+\intt_{t=2}^{x}{\frac{1}{t\left(\log{t}\right)^2}\summe_{p\leq t}{\e{-\alpha p}\log{p}}~\dd t}
\end{align*}
and apply afterwards the estimate
\begin{equation*}
\summe_{p\leq v}{\e{-\alpha p}\log{p}}=\frac{\mu(q)}{\varphi(q)}\summe_{n\leq v}{\e{-\beta n}}+{\rm O}\left(\frac{v}{\left(\log{v}\right)^A}\right)
\end{equation*}
that is valid for all $v,A\el\IR_{>1}$, see \cite[Lemma~3.1]{Vaughan}.

We then get
\begin{align*}
\T(\alpha)&=\frac{\mu(q)}{\varphi(q)}\left(\frac{1}{\log{x}}\summe_{2\leq n\leq x}{\e{-\beta n}}+\intt_{t=2}^{x}{\frac{1}{t\left(\log{t}\right)^2}\summe_{n\leq t}{\e{-\beta n}}\,\dd t}\right)+{\rm O}\left(\frac{x}{\left(\log{x}\right)^A}\right)\\
&=\Ts(\alpha)+{\rm O}\left(\frac{x}{\left(\log{x}\right)^A}\right)
\end{align*}
\end{Beweis}
\end{Lemma}

\begin{Folgerung}\label{S_Quadrat_auf_Major_Arcs}
We have
\begin{equation}\label{S_Quadrat_auf_Major_Arcs_Aussage}
\intt_{\mathfrak{M}}{\abs{\S(\alpha)}^2\dd\alpha}=x\zweinorm{f}^2+{\rm o}\left(\frac{x}{\log{x}}\right)
\end{equation}
\begin{Beweis}
Set
\begin{equation*}
\intt_{\mathfrak{M}}{\abs{S(\alpha)}^2\dd\alpha}=\Sigma_1+\Sigma_2+\Sigma_3+\Sigma_4
\end{equation*}
with 
\begin{align*}
\Sigma_1&=\intt_{\mathfrak{M}}{\abs{\Ss(\alpha)}^2\dd\alpha}=\summe_{q\leq Q}{\ssumme{a\leq q}{\abs{\G(q,a)}^2\intt_{\abs{\beta}\leq\frac{Q}{x}}{\abs{\summe_{n\leq x}{\e{\beta n}}}^2\dd\beta}}}\\
\Sigma_2&=\summe_{q\leq Q}{\ssumme{a\leq q}{\intt_{\abs{\beta}\leq\frac{Q}{x}}{\abs{\summe_{b\leq q}{\e{\frac{ab}{q}}\Xi\left(x;q,b;\beta\right)}}^2\dd\beta}}}\\
\Sigma_3&=\summe_{q\leq Q}{\ssumme{a\leq q}{\intt_{\abs{\beta}\leq\frac{Q}{x}}{\Ss\left(\frac{a}{q}+\beta\right)\overline{\summe_{b\leq q}{\e{\frac{ab}{q}}\Xi\left(x;q,b;\beta\right)}}\,\dd\beta}}}\\
\Sigma_4&=\summe_{q\leq Q}{\ssumme{a\leq q}{\intt_{\abs{\beta}\leq\frac{Q}{x}}{\overline{\Ss\left(\frac{a}{q}+\beta\right)}\summe_{b\leq q}{\e{\frac{ab}{q}}\Xi\left(x;q,b;\beta\right)}\,\dd\beta}}}
\end{align*}
To evaluate $\Sigma_1$ we complete the integration limits to $\left[-\frac{1}{2},\frac{1}{2}\right]$ and use properties of the exponential function. For the error term that occured, we apply (\ref{Gleichung:Exp_von_Beta_Abschaetzung}) and remark that $\norm{\beta}=\abs{\beta}$ for $\abs{\beta}\leq\frac{1}{2}$. The series converges with Parseval's identity to the limit $\zweinorm{f}^2$.
\begin{align*}
\Sigma_1&=\summe_{q\leq Q}{\ssumme{a\leq q}{\abs{\G(q,a)}^2}}\left(\floor{x}-\intt_{\frac{Q}{x}<\abs{\beta}\leq\frac{1}{2}}{\abs{\summe_{n\leq x}{\e{\beta n}}}^2\dd\beta}\right)\\
&=x\zweinorm{f}^2-x\summe_{q>Q}{\ssumme{a\leq q}{\abs{\G(q,a)}^2}}+{\rm O}\left(\frac{x}{Q}\right)
\end{align*}
The condition (\ref{Gleichung:Forderung_an_Betragsquadrat_ueber_Gausssche_Summe}) of theorem~\ref{Hauptsatz_der_DA} implies now together with the definition of $Q$
\begin{equation*}
\Sigma_1=x\zweinorm{f}^2+{\rm o}\left(\frac{x}{\log{x}}\right)
\end{equation*}
For $\Sigma_2$ we note that the function $\Xi$ defined in (\ref{Def_von_Xi}) satisfies the following inequality:
\begin{equation}\label{Gleichung:Abschaetzung_fuer_Xi}
\Xi(x;q,b;\beta)\ll\left(1+\abs{\beta}x\right)\max_{1\leq k\leq x}{\abs{\E(k;q,b)}}
\end{equation}
By neglecting the condition on co-primality for the sum over $a$, we get
\begin{align*}
\Sigma_2&\leq\summe_{q\leq Q}{\summe_{b_1,b_2\leq q}{~\intt_{\abs{\beta}\leq\frac{Q}{x}}{\Xi(x;q,b_1;\beta)\,\overline{\Xi(x;q,b_2;\beta)}\,\summe_{a\leq q}{\e{\frac{a}{q}\left(b_1-b_2\right)}}\dd\beta}}}\\
&=\summe_{q\leq Q}{q\summe_{b\leq q}{~\intt_{\abs{\beta}\leq\frac{Q}{x}}{\abs{\Xi(x;q,b;\beta)}^2\dd\beta}}}\\
&\ll Q\intt_{\abs{\beta}\leq\frac{Q}{x}}{\left(1+\abs{\beta}x\right)^2\dd\beta}\max_{1\leq k\leq x}{\summe_{q\leq Q}{\summe_{b\leq q}{\abs{\E(k;q,b)}^2}}}
\end{align*}
As $f$ fulfills condition (\ref{ZWEITE_ANFORDERUNG}), we can use the Barban-Davenport-Halberstam statement for limit-periodic functions, theorem \ref{Satz:BDH_fuer_grenzperiodische_Funktionen}, and get for all $A\el\IR$
\begin{equation*}
\Sigma_2\ll\frac{Q^4}{x}\max_{1\leq k\leq x}{\left(Qk+\frac{k^2}{\left(\log{k}\right)^A}\right)}\ll Q^5+Q^4\frac{x}{\left(\log{x}\right)^A}
\end{equation*}
Therefore
\begin{equation*}
\Sigma_2={\rm o}\left(\frac{x}{\log{x}}\right)
\end{equation*}
Using the estimates above for $\Sigma_1$ and $\Sigma_2$ and applying the Cauchy-Schwarz-inequality, we have for all $A\el\IR$
\begin{align*}
\Sigma_3+\Sigma_4&\ll\summe_{q\leq Q}{\ssumme{a\leq q}{\intt_{\abs{\beta}\leq\frac{Q}{x}}{\abs{\Ss\left(\frac{a}{q}+\beta\right)\summe_{b\leq q}{\e{\frac{ab}{q}}\Xi\left(x;q,b;\beta\right)}}\,\dd\beta}}}\\
&\ll\left(\Sigma_1\right)^{\frac{1}{2}}\left(\Sigma_2\right)^{\frac{1}{2}}\ll x^{\frac{1}{2}}\left(\frac{x}{\left(\log{x}\right)^A}\right)^{\frac{1}{2}}=\frac{x}{\left(\log{x}\right)^{\frac{A}{2}}}
\end{align*}
So
\begin{equation*}
\Sigma_3+\Sigma_4={\rm o}\left(\frac{x}{\log{x}}\right)
\end{equation*}
\end{Beweis}
\end{Folgerung}

\begin{Lemma}\label{Lemma:Major_Arcs_Abschaetzung}
We have for all $A\el\IR$
\begin{equation}\label{Major_Arcs_Abschaetzung_1}
\intt_{\mathfrak{M}}{\S(\alpha)\left(\T(\alpha)-\Ts(\alpha)\right)\dd\alpha}\ll Q^{\frac{3}{2}}\frac{x}{\left(\log{x}\right)^A}
\end{equation}
and
\begin{equation}\label{Major_Arcs_Abschaetzung_2}
\intt_{\mathfrak{M}}{\Ts(\alpha)\left(\S(\alpha)-\Ss(\alpha)\right)\dd\alpha}\ll Q\frac{x}{\left(\log{x}\right)^A}
\end{equation}

\begin{Beweis}
The Cauchy-Schwarz-inequality can be applied to equation (\ref{Major_Arcs_Abschaetzung_1}) and we get
\begin{equation*}
\abs{\intt_{\mathfrak{M}}{\S(\alpha)\left(\T(\alpha)-\Ts(\alpha)\right)}\dd\alpha}\leq\left(\intt_{\mathfrak{M}}{\abs{\S(\alpha)}^2\dd\alpha}\intt_{\mathfrak{M}}{\abs{\T(\alpha)-\Ts(\alpha)}^2\dd\alpha}\right)^{\frac{1}{2}}
\end{equation*}
For the first factor we use a trivial estimate from equation~(\ref{S_Quadrat_auf_Major_Arcs_Aussage})
\begin{equation*}
\intt_{\mathfrak{M}}{\abs{\S(\alpha)}^2\dd\alpha}\ll x
\end{equation*}
For the second factor, we use lemma~\ref{Lemma:S_Gleich_S_Stern_Plus_Rest} and get for all $A\el\IR$
\begin{equation*}
\intt_{\mathfrak{M}}{\abs{\T(\alpha)-\Ts(\alpha)}^2\dd\alpha}\ll\summe_{q\leq Q}{\ssumme{a\leq q}{\intt_{\abs{\beta}\leq\frac{Q}{x}}{\frac{x^2}{\left(\log{x}\right)^{2A}}\,\dd\beta}}}\ll Q^3\frac{x}{\left(\log{x}\right)^{2A}}
\end{equation*}
which proves (\ref{Major_Arcs_Abschaetzung_1}).

For the second statement we use the approximation property from lemma~\ref{Lemma:S_Gleich_S_Stern_Plus_Rest}. The left-hand side in (\ref{Major_Arcs_Abschaetzung_2}) is then equal to 
\begin{align*}
&\summe_{q\leq Q}{\ssumme{a\leq q}{\intt_{\abs{\beta}\leq\frac{Q}{x}}{\frac{\mu(q)}{\varphi(q)}\summe_{2\leq n\leq x}{\frac{\e{-\beta n}}{\log{n}}}\,\summe_{b\leq q}{\e{\frac{ab}{q}}\Xi(x;q,b;\beta)}\,\dd\beta}}}\\
=&\summe_{q\leq Q}{\summe_{b\leq q}{{\bf c}_q(b)\intt_{\abs{\beta}\leq\frac{Q}{x}}{\frac{\mu(q)}{\varphi(q)}\summe_{2\leq n\leq x}{\frac{\e{-\beta n}}{\log{n}}}~\Xi(x;q,b;\beta)\,\dd\beta}}}\\
\ll&\summe_{q\leq Q}{\summe_{b\leq q}{\abs{{\bf c}_q(b)}\intt_{\abs{\beta}\leq\frac{Q}{x}}{\abs{\frac{\mu(q)}{\varphi(q)}\summe_{2\leq n\leq x}{\frac{\e{-\beta n}}{\log{n}}}}~\abs{\Xi(x;q,b;\beta)}\,\dd\beta}}}
\end{align*}

Using the approximation (\ref{Gleichung:Ts_Abschaetzung}) for $\Ts$ and (\ref{Gleichung:Abschaetzung_fuer_Xi}) for $\Xi$ we get
\begin{align*}
&\summe_{q\leq Q}{\summe_{b\leq q}{\abs{{\bf c}_q(b)}\intt_{\abs{\beta}\leq\frac{Q}{x}}{\abs{\frac{\mu(q)}{\varphi(q)}\summe_{2\leq n\leq x}{\frac{\e{-\beta n}}{\log{n}}}}~\abs{\Xi(x;q,b;\beta)}\,\dd\beta}}}\\
\ll&\summe_{q\leq Q}{\frac{\mu(q)^2}{\varphi(q)}\summe_{b\leq q}{\abs{{\bf c}_q(b)}\intt_{\abs{\beta}\leq\frac{Q}{x}}{\frac{x}{\log{x}}\frac{1}{\left(1+\abs{\beta}x\right)}~\left(1+\abs{\beta}x\right)\max_{1\leq k\leq x}{\abs{\E(k;q,b)}}\,\dd\beta}}}\\
\ll&\frac{Q}{\log{x}}\max_{1\leq k\leq x}{\summe_{q\leq Q}{\frac{\mu(q)^2}{\varphi(q)}\summe_{b\leq q}{\abs{{\bf c}_q(b)}\abs{\E(k;q,b)}}}}
\end{align*}
Exploiting standard properties of Ramanujan's sum and the divisor function $\tf(q)$ results in
\begin{equation*}
\frac{Q}{\log{x}}\max_{1\leq k\leq x}{\summe_{q\leq Q}{\frac{\mu(q)^2}{\varphi(q)}q\,\tf(q)\max_{1\leq b\leq q}{\abs{\E(k;q,b)}}}}\ll Q\max_{1\leq k\leq x}{\summe_{q\leq Q}{\max_{1\leq b\leq q}{\abs{\E(k;q,b)}}}}
\end{equation*}
Finally, with applying condition (\ref{ZWEITE_ANFORDERUNG}) we get the desired result (\ref{Major_Arcs_Abschaetzung_2}).
\end{Beweis}
\end{Lemma}

\subsection{The main term}
\begin{Lemma}\label{Lemma:Major_Arcs_Hauptterm}On the major arcs we have
\begin{equation*}
\intt_{\mathfrak{M}}{\Ss(\alpha)\Ts(\alpha)\,\dd\alpha}=c_f\,\frac{x}{\log{x}}+{\rm o}\left(\frac{x}{\log{x}}\right)
\end{equation*}
with the absolute convergent series
\begin{equation}\label{C_f}
c_f:=\summe_{q=1}^{\infty}{\frac{\mu(q)}{\varphi(q)}\ssumme{a\leq q}{\G(q,a)}}
\end{equation}
\begin{Beweis}
Lemma~\ref{Lemma:MuQuadrat_durch_Phi_Summe} and the requirement~(\ref{Gleichung:Forderung_an_Betragsquadrat_ueber_Gausssche_Summe}) imply the absolute convergence of the series, as we have for all $v,w\el\IR_{>0}$ and $v':=\frac{\log{v}}{\log{2}}$, $w':=\frac{\log{w}}{\log{2}}-1$
\begin{align*}
\summe_{v<q\leq w}{\abs{\frac{\mu(q)}{\varphi(q)}\ssumme{a\leq q}{\G(q,a)}}}&\leq\summe_{v'\leq k\leq w'}{\summe_{2^k<q\leq 2^{k+1}}{\ssumme{a\leq q}{\abs{\frac{\mu(q)}{\varphi(q)}\G(q,a)}}}}\\
&\leq\summe_{v'\leq k\leq w'}{\left(\summe_{2^k<q\leq 2^{k+1}}{\frac{\mu(q)^2}{\varphi(q)}}\right)^{\frac{1}{2}}\left(\summe_{2^k<q\leq 2^{k+1}}{\ssumme{a\leq q}{\abs{\G(q,a)}^2}}\right)^{\frac{1}{2}}}\\
&\leq\summe_{k\geq v'}{\left(\summe_{q\leq 2^{k+1}}{\frac{\mu(q)^2}{\varphi(q)}}\right)^{\frac{1}{2}}\left(\summe_{q>2^k}{\ssumme{a\leq q}{\abs{\G(q,a)}^2}}\right)^{\frac{1}{2}}}\\
&\ll\summe_{k\geq v'}{\left(\left(\log{2^{k}}\right)\left(2^{-\frac{k}{r}}\right)\right)^{\frac{1}{2}}}\ll\summe_{k\geq v'}{\left(2^{-\frac{1}{4r}}\right)^k}~\stackrel{v\ra\infty}{\longrightarrow}{\,0}
\end{align*}
The number $r\el\IR_{>1}$ exists as we require (\ref{Gleichung:Forderung_an_Betragsquadrat_ueber_Gausssche_Summe}) to be true and it can be seen easily that the implicit constants can be chosen independently of $v$ and $w$. Cauchy's criterion implies the stated convergence.

To evaluate the integral
\begin{equation*}
\intt_{\mathfrak{M}}{\Ss(\alpha)\Ts(\alpha)\,\dd\alpha}=\summe_{q\leq Q}{\frac{\mu(q)}{\varphi(q)}\ssumme{a\leq q}{\G(q,a)}\intt_{\abs{\beta}\leq\frac{Q}{x}}{\summe_{m\leq x}{\e{\beta m}}\summe_{2\leq n\leq x}{\frac{\e{-\beta n}}{\log{n}}}\,\dd\beta}}
\end{equation*}
we complete the integration limits to $\left[-\frac{1}{2},\frac{1}{2}\right]$ and get
\begin{equation*}
\intt_{\mathfrak{M}}{\Ss(\alpha)\Ts(\alpha)\,\dd\alpha}=\frac{x}{\log{x}}\summe_{q\leq Q}{\frac{\mu(q)}{\varphi(q)}\ssumme{a\leq q}{\G(q,a)}}+\Sigma_5
\end{equation*}
with
\begin{equation*}
\Sigma_5\ll\abs{\summe_{q\leq Q}{\frac{\mu(q)}{\varphi(q)}\ssumme{a\leq q}{\G(q,a)}}}\left(\frac{x}{\left(\log{x}\right)^2}+\intt_{\frac{Q}{x}<\abs{\beta}\leq\frac{1}{2}}{\abs{\summe_{m\leq x}{\e{\beta m}}}\abs{\summe_{2\leq n\leq x}{\frac{\e{-\beta n}}{\log{n}}}}\dd\beta}\right)
\end{equation*}
With the convergence of the series over $q$ and lemma \ref{Lemma:Ts_Abschaetzung} as well as with the approximation (\ref{Gleichung:Exp_von_Beta_Abschaetzung}), we get
\begin{align*}
\intt_{\mathfrak{M}}{\Ss(\alpha)\Ts(\alpha)\,\dd\alpha}&=c_f\,\frac{x}{\log{x}}+{\rm O}\left(\frac{x}{\log{x}}\left(\frac{1}{\log{x}}+\abs{\summe_{q>Q}{\frac{\mu(q)}{\varphi(q)}\ssumme{a\leq q}{\G(q,a)}}}+\frac{1}{Q}\right)\right)\\
&=c_f\,\frac{x}{\log{x}}+{\rm o}\left(\frac{x}{\log{x}}\right)
\end{align*}
\end{Beweis}
\end{Lemma}

\begin{Folgerung}\label{Folgerung:Beitrag_der_Major_Arcs} On the major arcs we have
\begin{equation}\label{Gleichung:Beitrag_der_Major_Arcs}
\intt_{\mathfrak{M}}{\S(\alpha)\T(\alpha)\,\dd\alpha}=c_f\,\frac{x}{\log{x}}+{\rm o}\left(\frac{x}{\log{x}}\right)
\end{equation}
\begin{Beweis}
Writing
\begin{equation*}
\S(\alpha)\T(\alpha)=\S(\alpha)\left(\T(\alpha)-\Ts(\alpha)\right)+\left(\S(\alpha)-\Ss(\alpha)\right)\Ts(\alpha)+\Ss(\alpha)\Ts(\alpha)
\end{equation*}
and approximating the terms, yields the stated result.
\end{Beweis}
\end{Folgerung}

\subsection{The minor arcs}
\begin{Lemma}\label{Lemma:Beitrag_der_Minor_Arcs} For the integral on the minor arcs, we have
\begin{equation*}
\intt_{\mathfrak{m}}{\S(\alpha)\T(\alpha)\dd\alpha}={\rm o}\left(\frac{x}{\log{x}}\right)
\end{equation*}
\begin{Beweis}
We get with the Cauchy-Schwarz-inequality 
\begin{equation*}
\abs{\intt_{\mathfrak{m}}{\S(\alpha)\T(\alpha)\,\dd\alpha}}\leq\left(\intt_{\mathfrak{m}}{\abs{\S(\alpha)}^2\dd\alpha}\right)^{\frac{1}{2}}\left(\intt_{\mathcal{U}}{\abs{\T(\alpha)}^2\dd\alpha}\right)^{\frac{1}{2}}
\end{equation*}
An application of the prime number theorem yields then
\begin{equation*}
\intt_{\mathcal{U}}{\abs{\T(\alpha)}^2\dd\alpha}=\summe_{p\leq x}{1}\ll\frac{x}{\log{x}}
\end{equation*}
We get with
\begin{equation*}
\intt_{\mathfrak{m}}{\abs{\S(\alpha)}^2\dd\alpha}=\intt_{\mathcal{U}}{\abs{\S(\alpha)}^2\dd\alpha}-\intt_{\mathfrak{M}}{\abs{\S(\alpha)}^2\dd\alpha}=\summe_{n\leq x}{\abs{f(n)}^2}-\intt_{\mathfrak{M}}{\abs{\S(\alpha)}^2\dd\alpha}
\end{equation*}
and, luckily, as of condition (\ref{ERSTE_ANFORDERUNG}), we get with corollary~\ref{S_Quadrat_auf_Major_Arcs}:
\begin{align*}
\intt_{\mathfrak{m}}{\abs{\S(\alpha)}^2\dd\alpha}={\rm o}\left(\frac{x}{\log{x}}\right)
\end{align*}
\end{Beweis}
\end{Lemma}
Putting altogether: With equation (\ref{Gleichung:Hauptaussage_Zerlegung}), corollary \ref{Folgerung:Beitrag_der_Major_Arcs} and lemma~\ref{Lemma:Beitrag_der_Minor_Arcs} we get the statement (\ref{Gleichung:Haupt_Aussage}) and this completes the prove of theorem \ref{Hauptsatz_der_DA}.
\clearpage
\section{An application to $k$-free numbers}\label{Kapitel:k_freie_Zahlen_Beispielanwendung}
In this section we give an application of theorem \ref{Hauptsatz_der_DA}. For this purpose, let $s\el\IN$ and $\alpha_1,\dots\alpha_s\el\IN_0$, $r_1,\dots,r_s\el\IN$ with $2\leq r_1\leq\dots\leq r_s$ be fixed.

\begin{Definition}\label{EA_Def} For $a,q\el\IN$ we set the value of $\Ea{d_1,\dots,d_s,q}$ to $1$ (resp. $0$) if the following system of congruences 
\begin{equation}
\begin{aligned}\label{Kongruenzensystem}
n&\equiv -\alpha_j~\left(d_j\right)~~(1\leq j\leq s)\\
n&\equiv a~\left(q\right)
\end{aligned}
\end{equation}
has a solution in $n$ (resp. has no solution).
\end{Definition}

We choose our function $f$ to be
\begin{equation*}
f(n):=\mu_{r_1}(n+\alpha_1)\,\mu_{r_2}(n+\alpha_2)\cdot\ldots\cdot\mu_{r_s}(n+\alpha_s)
\end{equation*}
and 
\begin{equation*}
\mathfrak{F}:=\left\{n\el\IN~:~f(n)=1\right\}
\end{equation*}
The function $f$ is limit-periodic as is shown when using the lemmas \ref{muekistgrenzperiodisch} and~\ref{Lemma:Grenzp_und_beschr_Fkt} and the comments \ref{Bemerkungen_Grenzp_von_Shifted_Functions}. It only takes values from the set $\left\{0,1\right\}$ and satisfies the requirements from theorem \ref{Hauptsatz_der_DA} as will be shown below. We will apply similar methods as Brüdern \textit{et al.} \cite{Bruedern2}, \cite{Bruedern3} and Mirsky \cite{Mirsky1}, \cite{Mirsky2}.

To exclude the trivial case, we assume further the choice of the parameter ${\alpha_1,\dots\alpha_s}$, ${r_1,\dots,r_s}$ in such a way, that $\mathfrak{F}\neq\emptyset$. The following theorem characterizes exactly this case.

\begin{Satz}[Mirsky]\label{Satz_Mirsky}
The set $\mathfrak{F}$ is non-empty if and only if for every prime $p$ there exists a natural number $n$ with $\nmodu{n}{-\alpha_i}{p^{r_i}}$ for $1\leq i\leq s$. In this case, the set $\mathfrak{F}$ even has a positive density \cite[theorem 6]{Mirsky2}.
\end{Satz}

\subsection{Proof of the requirements (\ref{ERSTE_ANFORDERUNG}) and (\ref{ZWEITE_ANFORDERUNG})}
\begin{Definition}\label{Definition_von_D}
We define $\operatorname{D}(p)$ and $\operatorname{D}^{*}(p)$ as the number of natural numbers $n\leq p^{r_s}$ that solve at least one of the congruences $\modu{n}{-\alpha_i}{p^{r_i}}$ $(1\leq i\leq s)$, whereas we demand for $\operatorname{D}^{*}(p)$ in addition the condition $\ggt{n}{p}=1$, i.e.,
\begin{equation*}
\operatorname{D}(p):=\summe_{\substack{n\leq p^{r_s}\\\existsorig i:\,\moduu{n}{-\alpha_i}{p^{r_i}}}}{1}~~~~~~~~~~\operatorname{D}^{*}(p):=\summe_{\substack{n\leq p^{r_s}\\\existsorig i:\,\moduu{n}{-\alpha_i}{p^{r_i}}}}{\iverson{p\teiltnicht n}}
\end{equation*}
We set
\begin{equation}\label{Def_von_Frak_D}
\mathfrak{D}:=\pro_{p}{\left(1-\frac{\operatorname{D}(p)}{p^{r_s}}\right)}
\end{equation}
The convergence of this product follows from $\operatorname{D}(p)<p^{r_s}$ for every $p$ which is being implied by $\mathfrak{F}\neq\emptyset$, see \cite{Mirsky2}, and
\begin{equation}\label{Abschaetzung_fuer_D_1}
\operatorname{D}(p)\leq\summe_{i\leq s}{\summe_{n\leq p^{r_s}}{\iverson{\modu{n}{-\alpha_i}{p^{r_i}}}}}=\summe_{i\leq s}{p^{r_s-r_i}}\ll p^{r_s-2}
\end{equation}
\end{Definition}

\subsection*{The mean-value of $f$}
\begin{Satz}[Mirsky] We have for all ${\epsilon>0}$ 
\begin{equation*}
\summe_{n\leq x}{f(n)}=\mathfrak{D}\,x+{\rm O}\left(x^{\frac{2}{r_1+1}+\epsilon}\right)
\end{equation*}
See \cite[theorem 5]{Mirsky2}.
\end{Satz}
As the function $f$ can only assume the values $0$ or $1$, we also have
\begin{equation*}
\summe_{n\leq x}{\abs{f(n)}^2}=\mathfrak{D}\,x+{\rm O}\left(x^{\frac{2}{r_1+1}+\epsilon}\right)
\end{equation*}
and $\operatorname{M}(f)=\zweinorm{f}^2=\mathfrak{D}$. Therewith the requirement (\ref{ERSTE_ANFORDERUNG}) for $f$ follows.

\begin{Definition}
Set $g(q,a)$ as
\begin{align*}
g(q,a):=&\,\summe_{d_1,\dots,d_s=1}^{\infty}{\mu(d_1)\cdot\dots\cdot\mu(d_s)\frac{\Ea{d_1^{r_1},\dots,d_s^{r_s},q}}{\left[d_1^{r_1};\dots;d_s^{r_s}\right]}\left(\left[d_1^{r_1};\dots;d_s^{r_s}\right];q\right)}\\
=&\,q\summe_{d_1,\dots,d_s=1}^{\infty}{\mu(d_1)\cdot\dots\cdot\mu(d_s)\frac{\Ea{d_1^{r_1},\dots,d_s^{r_s},q}}{\left[d_1^{r_1};\dots;d_s^{r_s};q\right]}}
\end{align*}
It can be seen easily that the series converge.
\end{Definition}
\begin{Satz}
For $a,q\el\IN$ and $\epsilon>0$ we have
\begin{equation}\label{Gleichung:Restglied_fuer_Eta_f}
\summe_{\substack{n\leq x\\\moduu{n}{a}{q}}}{f(n)}=\frac{x}{q}g(q,a)+{\rm O}\left(x^{\frac{2}{r_1+1}+\epsilon}\right)
\end{equation}
whereas the implicit constant can be chosen independently from $a$ or $q$. The proof works analogous to the one in \cite[theorem 5]{Mirsky2}. It uses the identity (\ref{MueKIdentitaet}) for $\mu_k$ and 
\begin{equation*}
\summe_{\substack{n\leq x\\\moduu{n}{a}{q}}}{f(n)}=\summe_{\substack{n\leq x\\\moduu{n}{a}{q}}}{\summe_{d_1^{r_1}|(n+\alpha_1)}{\mu(d_1)}\dots\summe_{d_s^{r_s}|(n+\alpha_s)}{\mu(d_s)}}=\summe_{\substack{n\leq x\\\moduu{n}{a}{q}\\\forallorig j:\,\moduu{n}{-\alpha_j}{d_j^{r_j}}}}{\mu(d_1)\cdot\ldots\cdot\mu(d_s)}
\end{equation*}
\end{Satz}
Hence, the mean-value in residue classes is equal to ${\frac{1}{q}\,g(q,a)}$ and with the error term in (\ref{Gleichung:Restglied_fuer_Eta_f}) the validity of (\ref{ZWEITE_ANFORDERUNG}) for $f$ is shown.

This should be compared with the results of Brüdern~\textit{et al.} \cite{Bruedern3} and Brüdern \cite{Bruedern_AdditiveProblems}.

\subsection{The remainder of the singular series of $f$}
The validity of condition (\ref{Gleichung:Forderung_an_Betragsquadrat_ueber_Gausssche_Summe}) for $f$ is still open and will be shown in the following.
We start this section with an investigation of the function $g(q,a)$.

If we write
\begin{equation*}
g(q,a)=\summe_{d_1,\dots,d_s=1}^{\infty}{\theta_{a,q}(d_1,\dots,d_s)}
\end{equation*}
with
\begin{equation*}
\theta_{a,q}(d_1,\dots,d_s):=\mu(d_1)\cdot\ldots\cdot\mu(d_s)~\frac{\left(\left[d_1^{r_1};\dots;d_s^{r_s}\right];q\right)}{\left[d_1^{r_1};\dots;d_s^{r_s}\right]}~\Ea{d_1^{r_1},\dots,d_s^{r_s},q}
\end{equation*}
then $\theta_{a,q}(d_1,\dots,d_s)$ is a multiplicative function in $d_1,\dots,d_s$, which follows from the multiplicativity of the three factors
\begin{equation*}
\mu(d_1)\cdot\dots\cdot\mu(d_s),~~\frac{\left(\left[d_1^{r_1};\dots;d_s^{r_s}\right];q\right)}{\left[d_1^{r_1};\dots;d_s^{r_s}\right]},~~\Ea{d_1^{r_1},\dots,d_s^{r_s},q}
\end{equation*}

We then have
\begin{equation*}
g(q,a)=\pro_{p}{\chi_a^{(q)}(p)}~~~~~~~~~~\chi_a^{(q)}(p):=\summe_{\delta_1,\dots,\delta_s=0}^{\infty}{\theta_{a,q}(p^{\delta_1},\dots,p^{\delta_s})}
\end{equation*}
As $\mu(p^k)=0$ for $k\geq 2$ it follows
\begin{equation}\label{Gleichung_fuer_Chi}
\chi_a^{(q)}(p)=\summe_{\delta_1,\dots,\delta_s\el\left\{0,1\right\}}{\left(-1\right)^{\delta_1+\dots+\delta_s}\frac{\left(\left[p^{\delta_1r_1};\dots;p^{\delta_sr_s}\right];q\right)}{\left[p^{\delta_1r_1};\dots;p^{\delta_sr_s}\right]}~\Ea{p^{\delta_1r_1},\dots,p^{\delta_sr_s},q}}
\end{equation}
In the case $p\teiltnicht q$ we can write 
\begin{equation*}
\chi_a^{(q)}(p)=\left(1-\frac{\operatorname{D}(p)}{p^{r_s}}\right)
\end{equation*}
with \cite[theorem 5]{Mirsky2}. If we set in addition
\begin{equation}\label{Def_von_z_und_h}
z(q):=\pro_{p|q}{\left(1-\frac{\operatorname{D}(p)}{p^{r_s}}\right)^{-1}}~~~~~h(q,a):=\pro_{p|q}{\chi_a^{(q)}(p)}
\end{equation}
we get
\begin{equation*}
g(q,a)=\mathfrak{D}\,z(q)\,h(q,a)
\end{equation*}
and
\begin{equation}\label{Abschaetzung_fuer_z}
z(q)\ll 1
\end{equation}
by the comments in definition \ref{Definition_von_D}.

For $\modu{a}{b}{q}$ we have $h{(q,a)=h(q,b)}$.

\begin{Lemma}[Quasi-multiplicativity of $h$]
The function $h(q,a)$ is quasi-multiplicative, which means for all $q_1,q_2\el\IN$, $\left(q_1;q_2\right)=1$ and all $a_1,a_2\el\IN$ we have
\begin{equation*}
h(q_1q_2,a_1q_2+a_2q_1)=h(q_1,a_1q_2)\,h(q_2,a_2q_1)
\end{equation*}
The proof follows by elementary divisor relations.
\end{Lemma}

\begin{Definition}\label{Definition_von_H}
We set 
\begin{align*}
H(q,a)&:=\summe_{b\leq q}{h(q,b)\,\e{\frac{ab}{q}}}\\
\operatorname{H}(q)&:=\ssumme{a\leq q}{\abs{H(q,a)}^2}
\end{align*}
With the Gaußian sum $\G(q,a)=\summe_{b\leq q}{\frac{1}{q}g(q,b)\,\e{\frac{ab}{q}}}$ of $f$ it follows
\begin{equation}
\label{Gausssche_Summe_fuer_unser_f}
\ssumme{a\leq q}{\abs{\G(q,a)}^2}=\mathfrak{D}^2q^{-2}z(q)^2\H(q)
\end{equation}\end{Definition}

\subsection*{Properties of the function $\H$}
\begin{Lemma}\label{Lemma_Eigenschaften_von_H}
The function $\H$ has the following useful properties:
\begin{enumerate}
\item{\label{Item:Lemma_Aussage_1}$\H(q)$ is a multiplicative function}
\item{\label{Item:Lemma_Aussage_2}On prime powers we have
\begin{equation*}
\H(p^l)=\begin{cases}1&\mbox{for }l=0\\p^{3l-2r_s}\summe_{\substack{n,m\leq p^{r_s}\\\moduu{n}{m}{p^{l-1}}\\\existsorig i:\,\moduu{n}{-\alpha_i}{p^{r_i}}\\\existsorig i:\,\moduu{m}{-\alpha_i}{p^{r_i}}}}{\left(\iverson{\moduu{n}{m}{p^l}}-\frac{1}{p}\right)}&\mbox{for }1\leq l\leq r_s\\0&\mbox{for }l>r_s\\\end{cases}
\end{equation*}}
\item{\label{Item:Lemma_Aussage_3}We have the inequalities
\begin{equation*}
0\leq\H(p^l)\leq\begin{cases}s^2p^{3l-2r_1}&\mbox{for }1\leq l\leq r_1\\s^2p^{2l-r_1}&\mbox{for }r_1<l\leq r_s\end{cases}
\end{equation*}}
\end{enumerate}
\end{Lemma}
\begin{Beweis}
To statement \ref{Item:Lemma_Aussage_1}. Let $q_1,q_2\el\IN$, $\left(q_1;q_2\right)=1$ be given. We then have
\begin{align*}
\H(q_1q_2)&=\summe_{\substack{a\leq q_1q_2\\\left(a;q_1\right)=1\\\left(a;q_2\right)=1}}{\abs{H(q,a)}^2}=\ssumme{a_1\leq q_1}{\ssumme{a_2\leq q_2}{\abs{H(q_1q_2,a_1q_2+a_2q_1)}^2}}\\
&=\ssumme{a_1\leq q_1}{\ssumme{a_2\leq q_2}{\abs{\summe_{b\leq q_1q_2}{h(q_1q_2,b)\,\e{\frac{a_1q_2+a_2q_1}{q_1q_2}b}}}^2}}
\end{align*}
and by using the quasi-multiplicative property of $h$
\begin{equation*}
\begin{aligned}
\H(q_1q_2)&=\ssumme{a_1\leq q_1}{\ssumme{a_2\leq q_2}{\abs{\summe_{b_1\leq q_1}{\summe_{b_2\leq q_2}{h(q_1q_2,b_1q_2+b_2q_1)\,\e{\frac{a_1}{q_1}b_1q_2}\e{\frac{a_2}{q_2}b_2q_1}}}}^2}}\\
&=\ssumme{a_1\leq q_1}{\ssumme{a_2\leq q_2}{\abs{\summe_{b_1\leq q_1}{h(q_1,b_1q_2)\,\e{\frac{a_1}{q_1}b_1q_2}}}^2\abs{\summe_{b_2\leq q_2}{h(q_2,b_2q_1)\,\e{\frac{a_2}{q_2}b_2q_1}}}^2}}\\
&=\H(q_1)\H(q_2)
\end{aligned}
\end{equation*}
To statement \ref{Item:Lemma_Aussage_2}. We write
\begin{equation}
\begin{aligned}\label{Gleichung_H_hoch_pl}
\H(p^l)&=\summe_{a\leq p^l}{\abs{H(p^l,a)}^2}-\summe_{\substack{a\leq p^l\\p|a}}{\abs{H(p^l,a)}^2}\\
&=p^l\summe_{b\leq p^l}{h(p^l,b)^2}-p^{l-1}\summe_{\substack{b_1,b_2\leq p^l\\\moduu{b_1}{b_2}{p^{l-1}}}}{h(p^l,b_1)\,h(p^l,b_2)}
\end{aligned}
\end{equation}
and $\H(p^l)=0$ for $l>r_s$ follows, as in this case ${\modu{b_1}{b_2}{p^{r_s}}}$ and with the definition \ref{EA_Def} of $\EakeineKlammern$, the truth of $h(p^l,b_1)=h(p^l,b_2)$ is implied.
As $H$ is multiplicative, we have $\H(1)=1$. For the case $1\leq l\leq r_s$ we first evaluate $\chi_a^{(p^l)}(p)$:
\begin{equation*}
\begin{aligned}
\chi_a^{(p^l)}(p)&=p^l\summe_{\delta_1,\dots,\delta_s\left\{0,1\right\}}{\left(-1\right)^{\delta_1+\dots+\delta_s}\frac{\Ea{p^{\delta_1r_1},\dots,p^{\delta_sr_s},p^l}}{\left[p^{\delta_1r_1};\dots;p^{\delta_sr_s};p^l\right]}}\\
&=p^{l-r_s}\summe_{\substack{n\leq p^{r_s}\\\moduu{n}{a}{p^l}}}{\summe_{\substack{\delta_1,\dots,\delta_s\left\{0,1\right\}\\\forallorig j:\,\moduu{n}{-\alpha_i}{p^{\delta_ir_i}}}}{\left(-1\right)^{\delta_1+\dots+\delta_s}}}\\
&=p^{l-r_s}\summe_{\substack{n\leq p^{r_s}\\\moduu{n}{a}{p^l}}}{\summe_{\delta_1,\dots,\delta_s\left\{0,1\right\}}{\pro_{i\leq s}{\left(-1\right)^{\delta_i}\iverson{\modu{n}{-\alpha_i}{p^{\delta_ir_i}}}}}}\\
&=p^{l-r_s}\summe_{\substack{n\leq p^{r_s}\\\moduu{n}{a}{p^l}}}{\pro_{i\leq s}{\summe_{\delta_i\el\left\{0,1\right\}}{\left(-1\right)^{\delta_i}\iverson{\modu{n}{-\alpha_i}{p^{\delta_ir_i}}}}}}\\
\end{aligned}
\end{equation*}
The product has the value $0$ (resp. $1$) if $\modu{n}{-\alpha_i}{p^{r_i}}$ for at least one $i$ (resp. for no $i$ at all). Hence, we have for $1\leq l\leq r_s$
\begin{equation}\label{Chi_a_auf_kleinen_Primzahlpotenzen}
\chi_a^{(p^l)}(p)=p^{l-r_s}\summe_{\substack{n\leq p^{r_s}\\\moduu{n}{a}{p^l}\\\forallorig i:\,\nmoduu{n}{-\alpha_i}{p^{r_i}}}}{1}
\end{equation}
and with equation (\ref{Gleichung_H_hoch_pl}) and definition (\ref{Def_von_z_und_h}) of $h$
\begin{align}
\H(p^l)&=p^{3l-2r_s}\summe_{\substack{n,m\leq p^{r_s}\\\forallorig i:\,\nmoduu{n}{-\alpha_i}{p^{r_i}}\\\forallorig i:\,\nmoduu{m}{-\alpha_i}{p^{r_i}}}}{\summe_{\substack{b\leq p^l\\\moduu{b}{n}{p^l}\\\moduu{b}{m}{p^l}}}{1}~~-p^{3l-2r_s-1}\summe_{\substack{b_1\leq p^l\\\moduu{b_1}{n}{p^l}}}{\summe_{\substack{b_2\leq p^l\\\moduu{b_2}{m}{p^l}\\\moduu{b_1}{b_2}{p^{l-1}}}}{1}}}\nonumber\\
&=p^{3l-2r_s}\summe_{\substack{n,m\leq p^{r_s}\\\forallorig i:\,\nmoduu{n}{-\alpha_i}{p^{r_i}}\\\forallorig i:\,\nmoduu{m}{-\alpha_i}{p^{r_i}}}}{\iverson{\moduu{n}{m}{p^l}}}~\,-p^{3l-2r_s-1}\summe_{\substack{n,m\leq p^{r_s}\\\forallorig i:\,\nmoduu{n}{-\alpha_i}{p^{r_i}}\\\forallorig i:\,\nmoduu{m}{-\alpha_i}{p^{r_i}}}}{\iverson{\moduu{n}{m}{p^{l-1}}}}\nonumber\\
&=p^{3l-2r_s}\summe_{\substack{n,m\leq p^{r_s}\\\forallorig i:\,\nmoduu{n}{-\alpha_i}{p^{r_i}}\\\forallorig i:\,\nmoduu{m}{-\alpha_i}{p^{r_i}}}}{\left(\iverson{\moduu{n}{m}{p^l}}-\frac{1}{p}\iverson{\moduu{n}{m}{p^{l-1}}}\right)}\label{Gleichung_H_Zerlegung}
\end{align}
After multiple usages of
\begin{equation*}
\summe_{\substack{n\leq p^{r_s}\\\forallorig i:\,\nmoduu{n}{-\alpha_i}{p^{r_i}}}}{\iverson{\modu{n}{m}{p^u}}}=p^{r_s-u}-\summe_{\substack{n\leq p^{r_s}\\\existsorig i:\,\moduu{n}{-\alpha_i}{p^{r_i}}}}{\iverson{\modu{n}{m}{p^u}}}
\end{equation*}
we have for $0\leq u\leq r_s$:
\begin{equation*}
\summe_{\substack{n,m\leq p^{r_s}\\\forallorig i:\,\nmoduu{n}{-\alpha_i}{p^{r_i}}\\\forallorig i:\,\nmoduu{m}{-\alpha_i}{p^{r_i}}}}{\iverson{\modu{n}{m}{p^u}}}=p^{2r_s-u}-2p^{r_s-u}\summe_{\substack{n\leq p^{r_s}\\\existsorig i:\,\moduu{n}{-\alpha_i}{p^{r_i}}}}{1}+\summe_{\substack{n,m\leq p^{r_s}\\\existsorig i:\,\moduu{n}{-\alpha_i}{p^{r_i}}\\\existsorig i:\,\moduu{m}{-\alpha_i}{p^{r_i}}}}{\iverson{\modu{n}{m}{p^u}}}
\end{equation*}
If we specify $u=l$ and $u=(l-1)$ we get with equation (\ref{Gleichung_H_Zerlegung}):
\begin{equation}\label{Gleichung_Nummer2_H}
\H(p^l)=p^{3l-2r_s}\summe_{\substack{n,m\leq p^{r_s}\\\existsorig i:\,\moduu{n}{-\alpha_i}{p^{r_i}}\\\existsorig i:\,\moduu{m}{-\alpha_i}{p^{r_i}}}}{\left(\iverson{\modu{n}{m}{p^l}}-\frac{1}{p}\iverson{\modu{n}{m}{p^{l-1}}}\right)}
\end{equation}
The reader should compare the surprisingly similar representations of $\H(p^l)$ in (\ref{Gleichung_H_Zerlegung}) and (\ref{Gleichung_Nummer2_H}).

To statement \ref{Item:Lemma_Aussage_3}. With the definition~\ref{Definition_von_H} of $\H$ we always have $\H(p^l)\geq 0$ and with (\ref{Gleichung_Nummer2_H}) we have
\begin{align*}
\H(p^l)&\leq p^{3l-2r_s}\summe_{\substack{n,m\leq p^{r_s}\\\existsorig i:\,\moduu{n}{-\alpha_i}{p^{r_i}}\\\existsorig i:\,\moduu{m}{-\alpha_i}{p^{r_i}}}}{\iverson{\modu{n}{m}{p^l}}}\leq p^{3l-2r_s}\summe_{v,w\leq s}{\summe_{\substack{n,m\leq p^{r_s}\\\moduu{n}{-\alpha_v}{p^{r_v}}\\\moduu{m}{-\alpha_w}{p^{r_w}}}}{\iverson{\modu{n}{m}{p^l}}}}\\
&\leq p^{3l}\summe_{v,w\leq s}{p^{-\max{\left(l,r_v\right)}-r_w}}\leq p^{3l-\max{\left(l,r_1\right)}-r_1}\summe_{v,w\leq s}{1}\leq s^{2}\,p^{3l-\max{\left(l,r_1\right)}-r_1}
\end{align*}
which completes the proof.
\end{Beweis}

\begin{Folgerung}\label{Folgerung_fuer_Konvergenz}
For all ${\epsilon>0}$ and ${U\el\IR_{>0}}$, we have
\begin{equation*}
\summe_{U<q\leq 2U}{q^{-2}z(q)^2\H(q)}\ll U^{-\frac{r_1-1}{r_s}+\epsilon}
\end{equation*}
\end{Folgerung}
\begin{Beweis}
Using (\ref{Abschaetzung_fuer_z}) we get $z(q)^2\ll 1$. We start with
\begin{equation*}
\mathfrak{H}(U):=\summe_{U<q\leq 2U}{q^{-2}z(q)^2\H(q)}\ll U^{\frac{1}{r_1}-1}\summe_{q\leq 2U}{q^{-1-\frac{1}{r_1}}\H(q)}
\end{equation*}
The lemma~\ref{Lemma_Eigenschaften_von_H} shows that $\H(q)=0$, if $q$ is not $(r_s+1)$-free. It can be seen easily that every $(r_s+1)$-free number $q$ possesses a unique representation ${q=q_{1}q_{2}^{2}\cdot\ldots\cdot q_{r_s}^{r_s}}$  with pairwise co-prime and squarefree natural numbers $q_i$. Using lemma \ref{Lemma_Eigenschaften_von_H} again, we get
\begin{align*}
U^{\frac{1}{r_1}-1}\summe_{q\leq 2U}{q^{-1-\frac{1}{r_1}}\H(q)}&=U^{\frac{1}{r_1}-1}\summe_{q\leq 2U}{\summe_{q_{1}q_{2}^{2}\cdots q_{r_s}^{r_s}=q}{~\pro_{l\leq r_s}{q_l^{-l-\frac{l}{r_1}}\pro_{p|q_l}{\H(p^l)}}}}\\
&\ll U^{\frac{1}{r_1}-1+\epsilon}\summe_{q_{1}q_{2}^{2}\cdots q_{r_s}^{r_s}\leq 2U}{~\pro_{l\leq r_1}{q_l^{2l-\frac{l}{r_1}-2r_1}}\pro_{r_1<l\leq r_s}{q_l^{l-\frac{l}{r_1}-r_1}}}
\end{align*}

To simplify notations, we set ${\nu:=r_{s}-1}$ and~${\tau(l):=\frac{l}{r_s}\left(r_s-\frac{r_s}{r_1}-r_1+1\right)}$, and get
\begin{equation*}
\mathfrak{H}(U)\ll U^{\frac{1}{r_1}-1+\epsilon+\tau(1)}\summe_{q_{1}q_{2}^{2}\cdots q_{\nu}^{\nu}\leq 2U}{\pro_{l\leq r_1}{q_l^{2l-\frac{l}{r_1}-2r_1-\tau(l)}}\pro_{r_1<l\leq\nu}{q_l^{l-\frac{l}{r_1}-r_1-\tau(l)}}}
\end{equation*}
\begin{equation*}
\mathfrak{H}(U)\ll U^{-\frac{r_1-1}{r_s}+\epsilon}\summe_{q_{1}q_{2}^{2}\cdots q_{\nu}^{\nu}\leq 2U}{\pro_{l\leq r_1}{q_l^{l+l\frac{r_{1}-1}{r_s}-2r_1}}\pro_{r_1<l\leq\nu}{q_l^{l\frac{r_{1}-1}{r_s}-r_1}}}
\end{equation*}
As the exponents can't be larger than $-1$, the sums over $q_1,\dots,q_{\nu}$ are ${\rm O}(U^{\epsilon})$. If ${r_1<r_s}$, they are even convergent, which completes the proof.
\end{Beweis}

Now we are finally able to estimate the remainder of the singular series: The function $f$ satisfies the requirement (\ref{Gleichung:Forderung_an_Betragsquadrat_ueber_Gausssche_Summe}) with $r:=2r_s$ as with (\ref{Gausssche_Summe_fuer_unser_f}) and corollary \ref{Folgerung_fuer_Konvergenz}:
\begin{align*}
\summe_{q>\operatorname{w}}{\ssumme{a\leq q}{\left|\G_f(q,a)\right|^2}}&=\mathfrak{D}^2\summe_{q>\operatorname{w}}{q^{-2}z(q)^2\H(q)}=\mathfrak{D}^2\summe_{j=0}^{\infty}{\mathfrak{H}(2^{j}\operatorname{w})}\\
&\ll\operatorname{w}^{-\frac{r_1-1}{r_s}+\epsilon}\summe_{j=0}^{\infty}{\left(2^{-\frac{r_1-1}{r_s}+\epsilon}\right)^{j}}={\rm o}\left(\operatorname{w}^{-\frac{1}{2r_s}}\right)
\end{align*}
As all requirements of theorem \ref{Hauptsatz_der_DA} are fulfilled, we have
\begin{equation*}
\summe_{p\leq x}{f(p)}=c_f\,\frac{x}{\log{x}}+{\rm o}\left(\frac{x}{\log{x}}\right)
\end{equation*}
with the absolute convergent series
\begin{equation*}
c_f=\mathfrak{D}\summe_{q=1}^{\infty}{\frac{\mu(q)}{\varphi(q)}\,\frac{z(q)}{q}\,\ssumme{a\leq q}{H(q,a)}}
\end{equation*}

\subsection{Evaluation of the series $c_f$}
The function $q\mapsto\ssumme{a\leq q}{H(q,a)}$ is multiplicative. Let $q_1,q_2\el\IN$, $\left(q_1;q_2\right)=1$ be given. Then
\begin{align*}
\ssumme{a\leq q_1q_2}{H(q_1q_2,a)}&=\ssumme{a\leq q_1q_2}{\summe_{b\leq q_1q_2}{h(q_1q_2,b)\,\e{\frac{ab}{q_1q_2}}}}\\
&=\ssumme{a_1\leq q_1}{\ssumme{a_2\leq q_2}{\summe_{b_1\leq q_1}{\summe_{b_2\leq q_2}{h(q_1,b_1q_2)\,h(q_2,b_2q_1)\,\e{\frac{a_1}{q_1}b_1q_2}\,\e{\frac{a_2}{q_2}b_2q_1}}}}}\\
&=\left(\ssumme{a\leq q_1}{H(q_1,a)}\right)\left(\ssumme{a\leq q_2}{H(q_2,a)}\right)
\end{align*}

The other factors $\frac{\mu(q)}{\varphi(q)}$ and $\frac{z(q)}{q}$ in the representation of $c_f$ are trivially multiplicative. We then can write the series as an Euler product:
\begin{equation*}
c_f=\mathfrak{D}\pro_{p}{\left(\summe_{k=0}^{\infty}{\frac{\mu(p^k)}{\varphi(p^k)}\,\frac{z(p^k)}{p^k}\ssumme{a\leq p^k}{H(p^k,a)}}\right)}=\mathfrak{D}\pro_{p}{\left(1-\frac{z(p)}{p\left(p-1\right)}\ssumme{a\leq p}{H(p,a)}\right)}
\end{equation*}

Using properties of Ramanujan's sum, we have
\begin{equation}\label{Gleichung_fur_H_q_a}
\begin{aligned}
\ssumme{a\leq p}{H(p,a)}&=\summe_{b\leq p}{h(p,b)\ssumme{a\leq p}{\e{\frac{ab}{p}}}}=\summe_{b\leq p}{h(p,b)\,{\bf c}_p(b)}\\
&=\varphi(p)\,h(p,p)-\summe_{b<p}{h(p,b)}
\end{aligned}
\end{equation}
To evaluate ~${h(p,p)}$ and~${h(p,b)}$ we can apply the identity of (\ref{Chi_a_auf_kleinen_Primzahlpotenzen}) with $l=1$ and get
\begin{align*}
h(p,p)&=p^{1-r_s}\summe_{\substack{n\leq p^{r_s}\\\moduu{n}{b}{p}\\\forallorig i:\,\nmoduu{n}{-\alpha_i}{p^{r_i}}}}{1}=1-p^{1-r_s}\summe_{\substack{n\leq p^{r_s}\\\moduu{n}{0}{p}\\\existsorig i:\,\moduu{n}{-\alpha_i}{p^{r_i}}}}{1}
\end{align*}
as well as
\begin{align*}
\summe_{b<p}{h(p,b)}&=p^{1-r_s}\summe_{b<p}{\summe_{\substack{n\leq p^{r_s}\\\moduu{n}{b}{p}\\\forallorig i:\,\nmoduu{n}{-\alpha_i}{p^{r_i}}}}{1}}=\left(p-1\right)-p^{1-r_s}\summe_{b<p}{\summe_{\substack{n\leq p^{r_s}\\\moduu{n}{b}{p}\\\existsorig i:\,\moduu{n}{-\alpha_i}{p^{r_i}}}}{1}}\\
&=\left(p-1\right)-p^{1-r_s}\summe_{\substack{n\leq p^{r_s}\\\existsorig i:\,\moduu{n}{-\alpha_i}{p^{r_i}}}}{\summe_{\substack{b<p\\\moduu{b}{n}{p}}}{1}}\\
&=\left(p-1\right)-p^{1-r_s}\summe_{\substack{n\leq p^{r_s}\\\existsorig i:\,\moduu{n}{-\alpha_i}{p^{r_i}}}}{\left(1-\iverson{p|n}\right)}
\end{align*}

Using these results in (\ref{Gleichung_fur_H_q_a}) and noting that ${\varphi(p)=\left(p-1\right)}$, then
\begin{align*}
\ssumme{a\leq p}{H(p,a)}&=p^{1-r_s}\left(\operatorname{D}(p)-p\summe_{\substack{n\leq p^{r_s}\\\existsorig i:\,\moduu{n}{-\alpha_i}{p^{r_i}}}}{\iverson{p|n}}\right)
\end{align*}
and
\begin{equation*}
c_f=\mathfrak{D}\pro_{p}{\left(1-\frac{z(p)}{p^{r_s}\left(p-1\right)}\left(\operatorname{D}(p)-p\summe_{\substack{n\leq p^{r_s}\\\existsorig i:\,\moduu{n}{-\alpha_i}{p^{r_i}}}}{\iverson{p|n}}\right)\right)}
\end{equation*}

Looking again on the definitions (\ref{Def_von_Frak_D}) and~(\ref{Def_von_z_und_h}) of~$\mathfrak{D}$ and~$z(p)$, we get with those and ~${\varphi(p^{r_s})=p^{r_s-1}(p-1)}$
\begin{align*}
c_f&=\pro_{p}{\left(1-\frac{\operatorname{D}(p)}{p^{r_s}}-\frac{1}{p^{r_s}\left(p-1\right)}\left(\operatorname{D}(p)-p\summe_{\substack{n\leq p^{r_s}\\\existsorig i:\,\moduu{n}{-\alpha_i}{p^{r_i}}}}{\iverson{p|n}}\right)\right)}\\
&=\pro_{p}{\left(1-\frac{1}{\varphi(p^{r_s})}\left(\operatorname{D}(p)-\summe_{\substack{n\leq p^{r_s}\\\existsorig i:\,\moduu{n}{-\alpha_i}{p^{r_i}}}}{\iverson{p|n}}\right)\right)}\\
&=\pro_{p}{\left(1-\frac{\operatorname{D}^{*}(p)}{\varphi(p^{r_s})}\right)}
\end{align*}

The product is non-zero if and only if for each prime $p$ there exists a relatively prime natural number $n\leq p^{r_s}$ with ${\nmodu{n}{-\alpha_i}{p^{r_i}}}$ for all ${1\leq i\leq s}$, see theorem \ref{Satz_Mirsky}. In this case we have ${\operatorname{D}^{*}(p)<\varphi(p^{r_s})}$ and the convergence of the product is implied by ${\operatorname{D}^{*}(p)\ll p^{r_s-2}}$, see estimate (\ref{Abschaetzung_fuer_D_1}).

Hence, we have proven the identity (\ref{Gleichung:Aussage_des_Beispiels}).

\paragraph*{Acknowledgements.} The author wants to thank Jörg Brüdern who was the advisor of the author's diploma thesis created in 2007 on which this article is mainly based on.

\renewcommand{\refname}{References}
\bibliographystyle{abbrv}
\bibliography{Literatur}
\nocite{*}

\end{document}